\documentclass[11pt]{article}

\usepackage{times}
\usepackage{graphicx}
\usepackage{psfrag}
\usepackage{amsmath,amsthm}
\usepackage{amsfonts}

\theoremstyle{plain} %% This is the default
\newtheorem{thm}{Theorem}[section]

\newtheorem{lem}[thm]{Lemma}

\theoremstyle{definition}

\theoremstyle{remark}

\DeclareMathOperator*{\argmin}{arg\,min}
\DeclareMathOperator*{\argmax}{arg\,max}

\DeclareMathOperator{\Tr}{Tr}
\DeclareMathOperator{\re}{Re}
\DeclareMathOperator{\rank}{rank}

\begin{document}

\parskip=10pt

\flushbottom 

\title{Matrix product constraints by projection methods} 

\author{
Veit Elser\\
Department of Physics\\
Cornell University, Ithaca NY\\
} 

\date{}

\maketitle

\begin{abstract}
The decomposition of a matrix, as a product of factors with particular properties, is a much used tool in numerical analysis. Here we develop methods for decomposing a matrix $C$ into a product $X Y$, where the factors $X$ and $Y$ are required to minimize their distance from an arbitrary pair $X_0$ and $Y_0$. This type of decomposition, a projection to a matrix product constraint, in combination with projections that impose structural properties on $X$ and $Y$, forms the basis of a general method of decomposing a matrix into factors with specified properties. Results are presented for the application of these methods to a number of hard problems in exact factorization. 
\end{abstract}

\section{Introduction}

There is a large class of problems where the variables take the form of matrices $X$ and $Y$ that satisfy a product constraint
\begin{equation}\label{mpc}
X Y=C,
\end{equation}
as well as additional structural constraints that apply to $X$ and $Y$ individually. When the latter are ignored, and $X$ and $Y$ are completely unrestricted real or complex matrices, then it is easy, given $C$, to produce some decomposition of the form \eqref{mpc} where $X$ and $Y$ have a particular shape that is consistent with the shape and rank of $C$. However, such decompositions are far from unique and without additional properties are of little use in solving the complete problem, where the matrices must also satisfy structural constraints.

There is an additional, parameterized property we can impose on the decomposition \eqref{mpc} that will make it useful for solving the complete problem. This is the requirement that the decomposition minimizes the distances of the two matrices from an arbitrary pair $(X_0,Y_0)$. The decomposition is then said to be a projection of $(X_0,Y_0)$ to the matrix product constraint \eqref{mpc}. When combined with analogous projections that restore the structural constraints, the matrix product constraint projection makes available a variety of methods for solving the original problem. Although there are no solution guarantees when the problems are hard --- for which the constraint sets are non-convex -- projection methods as a heuristic are potentially useful because they can limit the search to matrices that are simultaneously close to both kinds of constraint.

\section{Simple projections for special factors}\label{sec:simple}

For three classes of factors the product constraint can be implemented directly on the original matrices $X$ and $Y$. The simple projections for these cases are discussed in this section. In the next section we will see how general product constraints can be reduced to constraints on combinations of special factors.

Our notation is appropriate for complex matrices but easily specializes to the real case by replacing the complex-conjugate transpose ($\dag$) with the transpose, unitary with orthogonal matrices etc. We let $U(m,n)$ denote unitary (orthogonal) matrices that are row unitary ($U U^\dag=I_m$) or column unitary ($U^\dag U=I_n$) for $m\le n$ or $m\ge n$, respectively.

\subsection{Symmetric factors}\label{sec:symmetric}

When $Y=X^\dag$ there is just one set of variables, $X\in \mathbb{C}^{m\times k}$ and
\begin{equation}
\|\Delta X\|^2_2=\Tr{(\Delta X^\dag \Delta X)}
\end{equation}
is the squared distance applied to the difference $\Delta X=X-X_0$ that the constraint projection minimizes. The constraint set is defined by
\begin{equation}
\mathcal{C}=\left\{X \in \mathbb{C}^{m\times k}\colon X X^\dag=C\right\},
\end{equation}
and the projection as
\begin{equation}\label{PC}
P_C(X_0)=\argmin_{X\in\mathcal{C}} \|X-X_0\|_2^2.
\end{equation}
Because the constraint set $\mathcal{C}$ is nonconvex, there will always be points $X_0$ for which there are multiple equally distant points $X$ on $\mathcal{C}$. In strict terms the projection is therefore a set-valued map. However, we will see that only for $X_0$ in a set  of measure zero is the distance minimizing point not unique.

An efficient method for computing \eqref{PC} has been known for a long time and arises, for example, when pairs of molecules (or models) are compared with allowance for arbitrary rotations to bring them into alignment.

To compute the projection $P_C$ of \eqref{PC} we obtain, as a one-time computation, the Cholesky decomposition of the constraint matrix $C=A A^\dag$, where $A\in\mathbb{C}^{m\times r}$ is lower triangular and $r=\rank{(C})\le \min{(m,k)}$. The constraint matrix $C$ is the Gram matrix of inner products of the rows of $X$, seen as vectors, and the rows of $A$ are a particular realization of $m$ vectors in a space of dimension $r$ that has the geometry implied by $C$.  The most general collection of $m$ vectors in a space of dimension $k\ge r$ that has the same geometry (Gram matrix) is given by
\begin{equation}\label{AU}
X=A U,
\end{equation}
where $U\in U(r,k)$. Computing the projection is thus an exercise in using the freedom in $U$ to minimize the distance between $X$ as defined by \eqref{AU} and an arbitrary matrix $X_0$.

By our definition of the squared distance the optimal $U$ is given by
\begin{eqnarray}
U&=&\argmin_{U'\in U(r,k)}\Tr{(A U'-X_0) (A U'-X_0)^\dag}\\
&=&\argmax_{U'\in U(r,k)}\re{\Tr{({X_0}^\dag A U')}}.
\end{eqnarray}
Expressing the singular value decomposition
\begin{equation}\label{VDW}
{X_0}^\dag A=V D W
\end{equation}
in terms of square unitary matrices $V\in U(k,k)$, $W\in U(r,r)$, the optimal $U$ is given by
\begin{eqnarray}
U&=&\argmax_{U'\in U(r,k)}\re{\Tr{(D W U' V)}}\\
&=&W^\dag \left(\argmax_{U''\in U(r,k)}\re{\Tr{(D U'')}} \right)V^\dag.\label{optU}
\end{eqnarray}
The diagonal matrix $D$ will have $r=\rank{(A)}$ positive values along the diagonal for a generic $X_0$, with $\rank{(X_0})\ge r$. When this is the case,
\begin{equation}
\re{\Tr{(D U'')}}=\sum_{i=1}^r D_{ii} \re{(U_{i i}'')},
\end{equation}
has a unique maximum among  $U''\in U(r,k)$ for $U_{i i}''=1$, $1\le i\le r$. Uniqueness is spoiled when $\rank{(X_0)}<r$, but this represents a set of measure zero. Comparing \eqref{optU} with \eqref{VDW}, we see that the projection can be compactly expressed as
\begin{equation}\label{symProj}
P_C(X_0)=A\, \mathcal{U}(A^\dag X_0),
\end{equation}
where the \textit{unitarization operator} $\mathcal{U}$ replaces all the singular values of a matrix by 1.

In the scalar case ($m=k=1$), where the constraint is $|x|^2=c$, the projection \eqref{symProj} reduces to
\begin{equation}\label{phaseProj}
P_c(x)=c \exp{(i\arg{x})}.
\end{equation}
This projection is used by almost all algorithms for solving the x-ray phase problem \cite{E1}.

Another simple case arises in searches for complex $m\times m$ Hadamard matrices \cite{TZ} $H$ defined by
\begin{equation}\label{Hadprod}
H H^\dag=m I_m,
\end{equation}
\begin{equation}\label{Hadstruc}
|H_{i j}|=1, \forall\, i,j.
\end{equation}
The projection to the product constraint \eqref{Hadprod} now simplifies to
\begin{equation}
P_C(H_0)=m\, \mathcal{U}(H_0),
\end{equation}
while the projection to the element-wise structure constraint \eqref{Hadstruc} is an instance of the scalar projection \eqref{phaseProj} with $c=1$.
For real Hadamard matrices the operator $\mathcal{U}$ acts on a real singular value decomposition (replacing all singular values by 1) and the structure projection is element-wise rounding to $\pm 1$.

\subsection{Orthogonal factors}\label{sec:orthogonal}

When $C=0$, the constraint $XY=0$ is geometrically the statement that the $m$ rows of $X$ and the $n$ columns of $Y$, seen as vectors, lie in orthogonal subspaces of $\mathbb{C}^k$. To project the pair $(X_0,Y_0)$ to this constraint set we must optimize both on the dimensions and the geometry of the orthogonal decomposition.

Let $r$, with $0\le r\le k$, be the dimension of the subspace into which the columns of $X_0$ are projected, then
\begin{equation}\label{XUU}
X=X_0 U U^\dag\qquad U\in U(k,r).
\end{equation}
The rows of $Y$ must then be in the subspace orthogonal to the one specified by $U$:
\begin{equation}\label{UUY}
Y=(I_k-U U^\dag)Y_0.
\end{equation}
Minimizing
\begin{equation}
\|X-X_0\|^2_2+\|Y-Y_0\|^2_2
\end{equation}
with respect to $r$ and $U$ defines the constraint projection $P_\perp(X_0,Y_0)=(X,Y)$. After some matrix manipulation, we arrive at the following:
\begin{equation}\label{Uopt}
U=\argmin_{U\in U(k,r),\;0\le r\le k}\Tr{\left((Y_0 {Y_0}^\dag-{X_0}^\dag X_0)U U^\dag\right)}.
\end{equation}
To solve the optimization problem we compute the eigen-decomposition
\begin{equation}
Y_0 {Y_0}^\dag-{X_0}^\dag X_0=V^\dag E V,
\end{equation}
where $V\in U(k,k)$ and $E$ is diagonal with real elements $E_{11}\le \cdots \le E_{kk}$. Since $V U=U'$ is again an arbitrary element of $U(k,r)$, we can rewrite \eqref{Uopt} as
\begin{equation}
U=V^\dag\left(\argmin_{U'\in U(k,r),\;0\le r\le k}\Tr{(E\,U' {U'}^\dag)}\right).
\end{equation}
Since the elements of $U' {U'}^\dag$ are always non-negative on the diagonal and bounded by $1$, the minimum is achieved when we select the first $r_-$ to be $1$ and the rest zero, where $r_-$ is the number of negative eigenvalues in $E$. The corresponding $U'$ will have $r_-$ columns and $1$'s on the diagonal, zero elsewhere. Relating this back to $U=V^\dag U'$ and \eqref{XUU}-\eqref{UUY}, we see that the projection can be written compactly as
\begin{eqnarray}
X&=&X_0\, \mathcal{E}_-(Y_0 {Y_0}^\dag-{X_0}^\dag X_0)\\
Y&=&\mathcal{E}_+(Y_0 {Y_0}^\dag-{X_0}^\dag X_0) Y_0,
\end{eqnarray}
where the \textit{eigenspace projection operators} $\mathcal{E}_\pm$ replace all the negative/positive eigenvalues by $1$, setting the rest to zero.

We are not aware of any applications that call for orthogonal matrix factors. However, we will see that the most general matrix product constraint (section \ref{sec:exrank}), when reduced to a form amenable by projections, calls for orthogonality in a decomposition of the factors as sums.

\subsection{Outer full rank factors}\label{sec:outerfullrank}

This is the core simple case upon which all (non-symmetric, $C\ne 0$) product constraint projections rely. To our knowledge the algorithm for this projection is new.

To be able to apply the simple projection derived in this section, the outer dimensions of the factors must match the rank of the constraint matrix: $m=n=\rank{(C)}=r$. In this section we therefore assume $C\in\mathbb{C}^{r\times r}$ is full rank and the factors have shapes $X\in\mathbb{C}^{r\times k}$, $Y\in\mathbb{C}^{k\times r}$, where $k\ge r$.
We wish to compute the projection
\begin{equation}
P_C(X_0,Y_0)=\argmin_{(X,Y)\in\mathcal{C}}\, \|X-X_0\|_2^2+\|Y-Y_0\|_2^2
\end{equation}
to the product constraint set
\begin{equation}
\mathcal{C}=\left\{(X,Y)\in \mathbb{C}^{r\times k}\times \mathbb{C}^{k\times r}\colon X Y=C\right\}.
\end{equation}

Our scheme for computing the projection is illustrated in Figure 1 for the simplest case of all: real matrices with $r=k=1$. While it is possible, in this scalar case, to obtain algebraic equations for the nearest point on the hyperbola, our method is iterative and generalizes to matrices. It comprises two operations: a \textit{quasiprojection} $Q$ and a true projection $P$ to the tangent-space approximation of the true constraint set.

\begin{figure}[t]
\begin{center}
\includegraphics[width=3.4in]{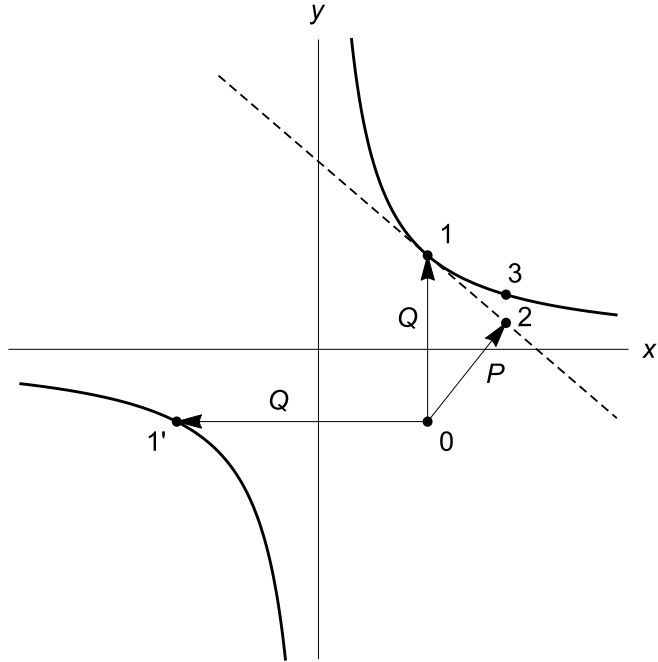}
\end{center}
\caption{Projection to the scalar product constraint $x y=c$ by iterating the quasiprojection $Q$ and the tangent space projection $P$. The quasiprojection takes the point to be projected, $(x_0,y_0)$, and constructs points $(x_1,y_1)$ and $(x_{1'},y_{1'})$ on the constraint set, selecting $(x_1,y_1)$ because it is closer to $(x_0,y_0)$. This is followed by $P$, which projects $(x_0,y_0)$ to the tangent space at $(x_1,y_1)$, producing point $(x_2,y_2)$. Another application of $Q$ produces the point $(x_3,y_3)$, an improvement over $(x_1,y_1)$ by its proximity to $(x_0,y_0)$.}
\end{figure}

The quasiprojection $Q(x_0,y_0;x,y)$ maps arbitrary pairs $(x,y)\in\mathbb{R}^2$ to the constraint set $x y = c$ by trying two alternatives and selecting the one that minimizes the distance to $(x_0,y_0)$. Starting with $(x,y)=(x_0,y_0)$, the two alternatives are $(c/y_0,y_0)$ and $(x_0,c/x_0)$. Whichever is closest to $(x_0,y_0)$ defines the first quasiprojection $(x_1,y_1)$. As is clear from Figure 1, $(x_1,y_1)$ is not the distance minimizing point on $x y=c$ to  $(x_0,y_0)$. To improve on $(x_1,y_1)$ we compute $P(x_0,y_0;x_1,y_1)=(x_2,y_2)$, a true distance minimizing point but on the tangent space approximation, at $(x_1,y_1)$, of the product constraint. This is followed by $Q(x_0,y_0;x_2,y_2)=(x_3,y_3)$ to bring the point back to the true constraint.

By iterating the two maps $T$ times, now for the general problem for complex matrices, we approximate the product constraint projection as
\begin{equation}
P_C(X_0,Y_0)\approx (X_T,Y_T),
\end{equation}
where
\begin{eqnarray}
(X_1,Y_1)&=&Q(X_0,Y_0;X_0,Y_0)\\
(X_{t+1},Y_{t+1})&=&Q(X_0,Y_0;P(X_0,Y_0;X_t,Y_t))),\quad 1\le t<T.\nonumber
\end{eqnarray}
For the intended applications of $P_C$, the point $(X_0,Y_0)$ maintains a respectful distance from the product constraint over most of the computation because of competing structural constraints. While it is important, via $Q$, to satisfy the product constraint precisely, the minimization of the (non-zero) distance brings diminishing returns. As we show later, in some problems even $T=2$ is adequate.

\subsubsection{Quasiprojection}\label{sec:quasiproj}

Given an arbitrary pair $(X_1,Y_1)$, our task here is to construct a pair $(X_2,Y_2)=Q(X_0,Y_0;X_1,Y_1)$ such that $X_2 Y_2=C$ and the distance between $(X_2,Y_2)$ and $(X_0,Y_0)$ is minimized when there are options. Our solution will have the other important property that when $(X_0,Y_0)\approx (X_1,Y_1)$ and $X_1 Y_1\approx C$, then $(X_2,Y_2)\approx (X_1,Y_1)$.

The two alternatives in the construction correspond to fixing $X_1$ or $Y_1$. Fixing $X_1$, we need to solve the equation
\begin{equation}
X_1 Y_2=C,
\end{equation}
for $Y_2$ or equivalently,
\begin{equation}\label{fixX}
X_1 \Delta Y=C-X_1 Y_0,
\end{equation}
for $\Delta Y=Y_2-Y_0$. Since $X_1\in \mathbb{C}^{r\times k}$ generically has full column rank, applying the Moore-Penrose pseudoinverse $X_1^+$ to \eqref{fixX} gives
\begin{equation}
\Delta Y=X_1^+(C-X_1 Y_0),
\end{equation}
the solution to \eqref{fixX} that minimizes $\|\Delta Y\|_2^2$. Therefore, the $X_1$-fixing option
\begin{equation}
(X_2,Y_2)=(X_1,Y_0+\Delta Y),
\end{equation}
gives squared distance
\begin{equation}\label{dist1}
\|X_1-X_0\|_2^2+\|\Delta Y\|_2^2.
\end{equation}
This is to be compared with fixing $Y_1$
\begin{eqnarray}
(X_2,Y_2)&=&(X_0+\Delta X,Y_1)\\
\Delta X&=&(C-X_0 Y_1)Y_1^+,
\end{eqnarray}
for which the squared distance is
\begin{equation}\label{dist2}
\|\Delta X\|_2^2+\|Y_1-Y_0\|_2^2.
\end{equation}
Whichever of \eqref{dist1} and \eqref{dist2} is smallest determines $(X_2,Y_2)$. The formulas for $\Delta X$ and $\Delta Y$ imply small changes, as required, when $(X_0,Y_0)\approx (X_1,Y_1)$ and both pairs approximately satisfy the product constraint.

\subsubsection{Tangent space projection}\label{sec:tangentproj}

For this projection we start with a pair $(X_1,Y_1)$ that satisfies $X_1 Y_1=C$ but is not necessarily distance minimizing to $(X_0,Y_0)$. The tangent space to the constraint at $(X_1,Y_1)$ is defined by pairs $(X,Y)$ that satisfy the linear equations
\begin{equation}\label{tanSpace}
(X-X_1)Y_1+X_1(Y-Y_1)=0.
\end{equation}
This is an $r\times r$ matrix of independent constraints, that we can impose on the problem of finding the distance minimizing point $(X_2,Y_2)$ using an $r\times r$ Lagrange multiplier matrix $F$:
\begin{equation}
(X_2,Y_2)_F=\argmin_{(X,Y)\in \mathbb{C}^{r\times k}\times \mathbb{C}^{k\times r}}\Tr G(X,Y;F)
\end{equation}
\begin{multline}
G(X,Y;F)=(X-X_0)^\dag(X-X_0)+(Y-Y_0)^\dag (Y-Y_0)\\
+F\left({Y_1}^\dag(X-X_1)^\dag+(Y-Y_1)^\dag {X_1}^\dag\right).
\end{multline}
Since $\Tr{G}$ is a positive definite quadratic form, it has a unique minimizer:
\begin{equation}\label{LMsol}
(X_2,Y_2)_F=(X_0-F {Y_1}^\dag,Y_0-{X_1}^\dag F).
\end{equation}

What remains is to find an $F$ such that $(X,Y)=(X_2,Y_2)_F$ satisfies the particular linear constraint \eqref{tanSpace}. Substituting \eqref{LMsol} into \eqref{tanSpace} we obtain the following equation for $F$:
\begin{equation}\label{Sylvester}
(X_1 {X_1}^\dag)F+F({Y_1}^\dag Y_1)=(X_0-X_1)Y_1+X_1(Y_0-Y_1).
\end{equation}
This is a Sylvester equation with positive definite coefficient matrices $A= X_1 {X_1}^\dag$ and $B={Y_1}^\dag Y_1$, since by satisfying  $X_1 Y_1=C$, $X_1$ and $Y_1$ both have rank $r$. For these conditions ($A$ and $B$ cannot have canceling eigenvalues) there is a unique solution for $F$ and therefore a unique distance minimizing projection \eqref{LMsol} to the tangent space constraint. There is a straightforward solution of the Sylvester equation for $F$ that starts with the singular value decompositions of $A$ and $B$.

\subsubsection{Product constraint projection for scalars}\label{subsubsec:scalar}

We record here, as a special case of the previous sections, the method for projecting to the product constraint for scalars. The formulas are given for a pair $(x,y)\in \mathbb{C}^2$ with constraint $x y=c\in \mathbb{C}$, but continue to hold  when these are real variables/constants. The complex conjugate of $x$ is written $\bar{x}$ and $x \bar{x}=|x|^2$.

To compute the quasiprojection $Q(x_0,y_0;x_1,y_1)=(x_2,y_2)$, we compare
\begin{equation}\label{xFix}
|x_1-x_0|^2+|c/x_1-y_0|^2
\end{equation}
with
\begin{equation}\label{yFix}
|c/y_1-x_0|^2+|y_1-y_0|^2.
\end{equation}
If \eqref{xFix} is less than \eqref{yFix}, $(x_2,y_2)=(x_1,c/x_1)$; otherwise, $(x_2,y_2)=(c/y_1,y_1)$.

To compute the tangent space projection $P(x_0,y_0;x_1,y_1)=(x_2,y_2)$ we note that the scalar case of the Sylvester equation \eqref{Sylvester} has the following solution for the scalar Lagrange multiplier $f$:
\begin{equation}
f=\frac{(x_0-x_1)y_1+(y_0-y_1)x_1}{|x_1|^2+|y_1|^2}.
\end{equation}
The projection to the tangent space is the scalar counterpart of \eqref{LMsol}:
\begin{equation}
(x_2,y_2)=(x_0-f \bar{y}_1,y_0-f \bar{x}_1).
\end{equation}

\section{Compound projections for general factors}\label{sec:compound}

The algorithms we use for decomposing $C$ into a product $XY$ where the factors also satisfy structural constraints require that all the constraints are implemented by just two projections. For the special types of factors in section \ref{sec:simple} this is done by imposing the product constraint, say for outer full rank factors, by the first projection,
\begin{equation}\label{P1simple}
P_1(X,Y)=P_C(X,Y)
\end{equation}
and all the structure constraints by the second projection:
\begin{equation}\label{P2simple}
P_2(X,Y)=(P_*(X),P_*(Y)).
\end{equation}
Here $P_*$ denotes the projection to a specific structure, and may be different for the two factors. In many applications the structure constraints are element-wise. For example, in non-negative matrix factorization we set $P_*=P_+$, the projection that sets all negative elements to zero and keeps the others unchanged.

In this section our goal is to again construct a pair of projections, such as \eqref{P1simple} and \eqref{P2simple}, but for factors not among the special types in section \ref{sec:simple}. We give three such constructions. The first two build on the projection for outer full rank factors and differ with respect to the ranks of the factors matching or exceeding the rank of $C$. Our third construction has no restrictions on the factors but is furthest in spirit from imposing a product constraint in that the product of $X$ and $Y$ is expressed as a sum of rank-1 matrices.

\subsection{Rank-limited factors}\label{sec:ranklim}

We now have $X\in \mathbb{C}^{m\times k}$ and $Y\in \mathbb{C}^{k\times n}$, and the knowledge that $\rank{(X)}=\rank{(Y)}=\rank{(C)}= r\le \min{(m,k,n)}$. If $m=n=r$ then the simple projection of section \ref{sec:outerfullrank} can be used and the construction described here is unnecessary. If just one of the outer dimensions matches $r$, the hybrid construction described at the end of this section should be used.

As a one-time computation we obtain the singular value decomposition of $C$,
\begin{equation}\label{Csvd}
C=U D V,
\end{equation}
where $U\in U(m,r)$, $V\in U(r,n)$, and $D$ is the diagonal matrix of the $r$ sorted singular values. To help tailor the projection method to specific applications, we introduce a two parameter rescaling in this decomposition: $U\to g\,U$, $V\to h\,V$, $D\to D/(g h)$. Henceforth we use the symbols $U$, $V$ and $D$ with this rescaling in effect, so that
\begin{equation}\label{UVscale}
U^\dag U=g^2 I_r\qquad V V^\dag =h^2 I_r.
\end{equation}

Since $X$ has rank $r$, the constraint $XY=C$ implies that $X$ is in the column-span of the $r$ columns of $U$. Similarly, $Y$ is in the row-span of $V$. We may therefore write
\begin{equation}\label{lin}
X=U W\qquad Y=Z V,
\end{equation}
where $W\in \mathbb{C}^{r\times k}$, $Z\in \mathbb{C}^{k\times r}$ satisfy the constraint
\begin{equation}\label{wz}
W Z = D.
\end{equation}
Given variable pairs $(W,Z)$ we can use the outer full rank projection of section \ref{sec:outerfullrank} to project to constraint \eqref{wz}.

To design projections that solve the original problem for the factors $X$ and $Y$ we work with the matrix pairs $W,X$ and $Z,Y$. Three kinds of constraints apply to these: the product constraint \eqref{wz}, the linear constraints \eqref{lin}, and structural constraints on $X$ and $Y$. A pair of compound projections that implements all of these constraints is the following,
\begin{eqnarray}
P_1(W,X;Z,Y)&=&(W',P_*(X);Z',P_*(Y))\label{P1ranklimited}\\
(W';Z')&=&P_C(W;Z)\\
P_2(W,X;Z,Y)&=&(P_U(W,X);P_V(Z,Y)),
\end{eqnarray}
where $P_C$ is the outer full rank projection with $C=D$, and $P_U$ and $P_V$ project to the linear constraints \eqref{lin}. To verify that this is a valid compound projection construction for the original problem we check two things. First, we note that in both $P_1$ and $P_2$ each of the variables appears at most once as the argument of a simple projection. The second check is to note that if $(W,X;Z,Y)$  is fixed by both $P_1$ and $P_2$ then (\textit{i}) $X$ and $Y$ have the correct structure, (\textit{ii}) have the correct product because the pairs $W,X$ and $Z,Y$ satisfy \eqref{lin} for a particular pair $(W,Z)$ that satisfies \eqref{wz}. We see that the singular value structure of the constraint matrix $C$ is exploited not just by the presence of the singular value matrix $D$ in the product projection $P_C$ (inside $P_1$), but also the corresponding column and row information in the projections $P_U$ and $P_V$ (inside $P_2$).

The projections to the linear compatibility constraints \eqref{lin}, though straightforward, bring up a question on the distance used in defining the projections. Because these operate in the Cartesian-product space comprising all four matrices, our choice of distance may want to respect intrinsic differences among them. In particular, when projecting to the constraint $X=U W$ one might want to define the squared distance by
\begin{equation}\label{gDist}
\|\Delta X\|_2^2 + g^2 \|\Delta W\|_2^2,
\end{equation}
with a freely adjustable metric parameter $g$. Alternatively, in terms of new matrices $U'= g U$ and  $W'= W/g$ the form of the linear constraint is unchanged but the parameter $g$ in the distance is eliminated. As this last option is more convenient, henceforth we use distances with an artificial symmetry among the different components ($g=1$ in \eqref{gDist}) and instead absorb the metric freedom in the definitions of $U$ and $V$. It is for this reason that we introduced the two-parameter rescaling of the standard singular value decomposition \eqref{Csvd} where $U$ and $V$ have normalizations \eqref{UVscale}.

To compute the projection $P_U(W_0,X_0)=(W_1,X_1)$, where $X_1=U W_1$, we only perform a minimization over $W$ since $X$ can be directly expressed in terms of $W$ when the constraint is satisfied:
\begin{equation}
W_1=\argmin_{W\in \mathbb{C}^{r\times k}}\Tr{\left((U W-X_0)^\dag (U W-X_0)+(W-W_0)^\dag(W-W_0)\right)}.
\end{equation}
Minimizing this positive definite quadratic form and using $U^\dag U=g^2 I_r$, we obtain
\begin{equation}
W_1=(W_0+U^\dag X_0)/(g^2+1), \qquad X_1=U W_1.
\end{equation}
Similarly,
\begin{equation}
Z_1=(Z_0+Y_0 V^\dag)/(h^2+1), \qquad Y_1=Z_1 V.
\end{equation}

In the event that one of the outer dimensions, say $n$, equals the rank $r$ we would use a simplified compound construction:
\begin{eqnarray}
P_1(W,X;Y)&=&(W',P_*(X);Y')\\
(W';Y')&=&P_C(W;Y)\\
P_2(W,X;Y)&=&(P_U(W,X);P_*(Y)).
\end{eqnarray}
The constraint matrix in projection $P_C$ is now $DV$, that is, only `half' of the singular value decomposition of the original constraint matrix. As in the general case, it is straightforward to verify the validity of this compound construction for solving the original problem.

\subsection{Rank-excessive factors}\label{sec:exrank}

This is the most elaborate case, but it does arise in applications. For example, the \textit{linear Euclidean distance matrix}
\begin{equation}
C=\left[
\begin{array}{cccccc}
0&1&4&9&16&25\\
1&0&1&4&9&16\\
4&1&0&1&4&9\\
9&4&1&0&1&4\\
16&9&4&1&0&1\\
25&16&9&4&1&0
\end{array}
\right]
\end{equation}
has rank $r=3$ and a non-negative factorization into $X\in\mathbb{R}^{6\times 5}$ and $Y\in\mathbb{R}^{5\times 6}$ \cite{GG}. The non-negative rank of $C$ is therefore bounded by $5$. However, the factors have excessive rank $4>r$ and therefore cannot be found with the compound construction of the previous section.

To treat this case we decompose the factors first as sums:
\begin{equation}
X=X_C+X_\perp\qquad Y=Y_C+Y_\perp.
\end{equation}
Here $X_C$ and $Y_C$ are to be interpreted as the parts of the factors that participate in the product while $X_\perp$ and $Y_\perp$ roughly correspond to what is left over. In more precise terms, we define $X_C$ and $Y_C$ exactly as we would in the full rank case:
\begin{equation}\label{linex}
X_C=U W\qquad Y_C=Z V.
\end{equation}
By construction, $X_C$ and $Y_C$ have rank $r=\rank{(C)}$ and product $X_C Y_C=UDV=C$ when $W Z=D$. The parts $X_\perp$ and $Y_\perp$ make up for the excess rank.

The original product constraint, $X Y=C$, implies the following constraint on the parts:
\begin{equation}\label{ortho3}
X_C Y_\perp+X_\perp Y_C+X_\perp Y_\perp =0.
\end{equation}
We can project to this constraint, in a compound setting, by introducing replicated variables \cite{GE} $\widetilde{X}_C$, $\widetilde{X}_\perp$, $\widetilde{Y}_C$ and $\widetilde{Y}_\perp$. As their name suggests, replicated variables satisfy the simple equality constraints:
\begin{equation}
X_C=\widetilde{X}_C\quad X_\perp=\widetilde{X}_\perp\quad Y_C=\widetilde{Y}_C\quad Y_\perp=\widetilde{Y}_\perp.
\end{equation}
In fact, these constraints are so simple that they can be combined with the projection to the structure constraints. We therefore write the structure projection in the expanded form
\begin{equation}
P_*(X_C,\widetilde{X}_C,X_\perp,\widetilde{X}_\perp)=(X'_C,X'_C,X'_\perp,X'_\perp)
\end{equation}
where $X'_C+X'_\perp$ satisfies the structural constraint on the original matrix $X$. Computing this projection for element-wise structure constraints is easy as it only involves four numbers at a time.

Non-negativity of $X$ would be treated in the following way. Suppose $x_C$, $\tilde{x}_C$, $x_\perp$ and $\tilde{x}_\perp$ are the four real scalar elements on which we want to compute the projection $P_*$. The first step is to project to the equality constraints $x'_C=\tilde{x}'_C=\bar{x}_C=(x_C+\tilde{x}_C)/2$ and $x'_\perp=\tilde{x}'_\perp=\bar{x}_\perp=(x_\perp+\tilde{x}_\perp)/2$. Now if $\bar{x}_C+\bar{x}_\perp>0$ we are done and the result of the projection is $(\bar{x}_C,\bar{x}_C,\bar{x}_\perp,\bar{x}_\perp)$. If that is not the case, we shift both parts by the same amount to give a sum of zero; the resulting projection is $(\delta x, \delta x, -\delta x,-\delta x)$, where $\delta x=(\bar{x}_C-\bar{x}_\perp)/2$.

Since the variables $W$ and $Z$ do not appear in the structure constraints, we combine them as in \eqref{P1ranklimited} when forming the first compound projection:
\begin{multline}
P_1(W,X_C,\widetilde{X}_C,X_\perp,\widetilde{X}_\perp;Z,Y_C,\widetilde{Y}_C,Y_\perp,\widetilde{Y}_\perp)=\\
(W',P_*(X_C,\widetilde{X}_C,X_\perp,\widetilde{X}_\perp); Z', P_*(Y_C,\widetilde{Y}_C,Y_\perp,\widetilde{Y}_\perp))\\
(W';Z')=P_C(W;Z).
\end{multline}
Having replicas of $X_C$, $X_\perp$, $Y_C$ and $Y_\perp$ makes it possible to project to the remaining constraints, \eqref{linex} and \eqref{ortho3}. These can be written in terms of replicas such that no variable appears in more than one constraint:
\begin{equation}
X_C=U W\qquad Y_C=Z V\label{con2a}
\end{equation}
\begin{equation}\label{con2b}
\widetilde{X}_C Y_\perp+X_\perp \widetilde{Y}_C+\widetilde{X}_\perp \widetilde{Y}_\perp =0.
\end{equation}
Projecting to constraint \eqref{con2a} is accomplished with the same projections $P_U$ and $P_V$ that are used in the rank-limited case. Constraint \eqref{con2b} is an instance of orthogonal factors (section \ref{sec:orthogonal}), as is clear when we column-concatenate $\widetilde{X}_C$, $X_\perp$ and $\widetilde{X}_\perp$ to form $X_3\in \mathbb{C}^{m\times 3k}$ and row-concatenate $Y_\perp$, $\widetilde{Y}_C$,  and $\widetilde{Y}_\perp$ to form $Y_3\in \mathbb{C}^{3k\times n}$ (for constraint $X_3 Y_3 =0$). The second compound projection is therefore
\begin{multline}
P_2(W,X_C,\widetilde{X}_C,X_\perp,\widetilde{X}_\perp;Z,Y_C,\widetilde{Y}_C,Y_\perp,\widetilde{Y}_\perp)=\\
(P_U(W,X_C),\widetilde{X}'_C,X'_\perp,\widetilde{X}'_\perp); P_V(Z, Y_C),\widetilde{Y}'_C,Y'_\perp,\widetilde{Y}'_\perp))\\
(\widetilde{X}'_C,X'_\perp,\widetilde{X}'_\perp;\widetilde{Y}'_C,Y'_\perp,\widetilde{Y}'_\perp)=P_\perp(\widetilde{X}_C,X_\perp,\widetilde{X}_\perp;\widetilde{Y}_C,Y_\perp,\widetilde{Y}_\perp).
\end{multline}
It is easy to check that if both $P_1$ and $P_2$ fix all ten matrix variables, then $X=X_C+X_\perp$ and $Y=Y_C+Y_\perp$ have the correct product and satisfy the structure constraints. There is an exchange of information between the two factors in both $P_1$ and $P_2$, while this is true only for $P_1$ in the rank-limited case (which operates on only four matrix variables).

\subsection{Rank-1 decomposition}\label{sec:rank1}

The matrix product constraint \eqref{mpc} can be written in the form
\begin{equation}
\sum_{l=1}^k Z^l=C,
\end{equation}
where the $Z^l\in \mathbb{C}^{m\times n}$ are required to be rank-1 matrices:
\begin{equation}
Z^l={x^l}^\dag \,y^l\qquad 1\le l\le k.
\end{equation} 
Here $x^l\in \mathbb{C}^m$, $y^l\in \mathbb{C}^n$ are row vectors. The difficulty of recovering $x^l$ and $y^l$ from $Z^l$ (the explicit factors $X$ and $Y$) may depend on the nature of the structure constraints. The non-negativity constraint represents an easy case. For suppose we have a solution of real, non-negative and rank-1 $Z$'s that sum to $C$. To decompose $Z^l$ as $Z^l _{i j}=x^l_i\, y^l_j$ into nonnegative vectors $x^l$ and $y^l$ we can proceed as follows. Find an $i$ for which the row $Z^l _{i j}$ is not entirely zero, set $x^l_i=a^l>0$ and thereby infer $y^l_j=Z^l _{i j}/a^l$ for all $j$. Now take a $j$ for which $y^l_j>0$ and determine $x^l_i=Z^l _{i j}/y^l_j$ for all $i$. In this way one obtains matrix factors $X$ and $Y$ with $k$ arbitrary scale parameters $a^l$ and a permutation arbitrariness in how the $k$ summands are assigned to the $k$ rows/columns of the factors.

In this projection scheme the variables are $Z^1,\ldots, Z^k$ and again there are two projections that act in this space. The first projection acts on the $Z$'s individually,
\begin{equation}\label{rank1P1}
P_1(Z^1,\ldots,Z^k)=(P_{r1}(Z^1),\ldots,P_{r1}(Z^k)),
\end{equation}
with $P_{r1}$ projecting each to the nearest rank-1 matrix. The second projection combines structural constraints with the constraint that the $Z$'s have sum $C$:
\begin{equation}\label{rank1P2}
P_2(Z^1,\ldots,Z^k)=P_*(Z^1,\ldots,Z^k).
\end{equation}
Because most structure constraints are element-wise, the computation of $P_*$ is usually only slightly more complicated than projecting to the structure constraints without the property that the sum is $C$. The case of non-negativity is worked out below.

The algorithm for computing $P_{r1}$ is well known and is concisely described as setting to zero all but the largest singular value of the matrix, everything else being left unchanged.
As non-negativity is a widely used structure constraint, we devote the rest of this section to the algorithm for computing the projection $P_*$ to this constraint. The constraints associated with a particular matrix element of $C$ have the form
\begin{equation}\label{simplex}
\sum_{l=1}^k z^l=c,\qquad z^l\ge 0,\;\forall\, l,
\end{equation}
where $z^l=Z^l_{i j}$ and $c=C_{i j}\ge 0$ are the variables and constant that go with the $(i,j)$ matrix element. We can refer to $P_*$ as the \textit{simplex projection} because the set of feasible $k$-tuples for \eqref{simplex} forms a regular $k-1$ simplex. To implement non-negativity, the simplex projection is applied independently on the $k$-tuples at each $(i,j)$.

An efficient computation of $P_*$ is based on two simple lemmas that we state without proof. This projection is built from two simpler projections that act on $k$-tuples $z$: $P_c(z)$ projects to the constraint \eqref{simplex} with non-negativity relaxed (all variables are shifted by the same value so as to produce the correct sum), and $P_0(z)$ replaces $z$ by all zeros. In our notation all three operators ($P_*$, $P_c$, $P_0$) continue to act on direct sums of arbitrary subsets the original variables, with no change in the value of $c$.
\begin{lem}
For $1\le k'\le k$ and all $z\in\mathbb{R}^{k'}$,
\begin{equation}
P_*(z)=P_*(P_c(z)).
\end{equation}
\end{lem}
\begin{lem}
For $1\le k'\le k$ and all $z\in\mathbb{R}^{k'}$, let $P_c(z)=z_{c+}\oplus z_{c-}$ be the direct sum decomposition into positive and nonpositive values; then
\begin{equation}
P_*(z_{c+}\oplus z_{c-})=P_*(z_{c+})\oplus P_0(z_{c-}).
\end{equation}
\end{lem}
\noindent In combination, the two lemmas give an efficient recursive algorithm for $P_*$. To efficiently manage the direct sums (partitioning into positive and nonpositive values) the initial $z$ should first be permuted into a sorted order.

When $C$ is an integer matrix and we believe there is a rank-1 decomposition where all the $Z$'s are also integer matrices, we can use the stronger structure constraint where all the $z$'s in \eqref{simplex} are required to be non-negative integers. To project to this constraint we use the composition $P_A\circ P_*$, where $P_*$ is the simplex projection for sum $c$ as above, and $P_A$ is the projection \cite{CS} to the $A_{k-1}$ root lattice (suitably shifted so the $k$-tuples sum to $c$ rather than zero). Establishing that this is a projection requires a check that the simplex of the first projection is covered by lattice Voronoi cells belonging only to lattice points that lie in the simplex.

While the rank-1 method comes without restrictions on the factors, and the constraint projections are relatively easy to compute, there are two reasons to favor the alternative method that uses the projection $P_C$. First, the rank-1 method treats the constraint matrix $C$ as a structureless set of $m n$ numbers. By contrast, the methods in sections \ref{sec:ranklim} and \ref{sec:exrank} exploit the singular value structure of $C$ which surely is advantageous when $C$ is dominated by a few singular values. Second, the rank-1 method requires significantly more variables: $m n k$ compared to $(m+n+2r)k$ (rank-limited) or $(4m+4n+2r)k$ (rank-excessive).

\section{Constraint satisfaction by iterated projections}\label{sec:iterproj}

In all the projection methods described above, simple or compound, the variables are Cartesian products of various complex or real matrices. For the purposes of this section we can treat these as vectors $x\in \mathbb{C}^M$ or $x\in \mathbb{R}^M$, where $M$ is the total number of variables in the Cartesian product. Also, solutions $x^*$ to all problems are identified by the property that they are fixed by two projections:
\begin{equation}
P_1(x^*)=x^*\qquad P_2(x^*)=x^*.
\end{equation}
The convention of the preceding sections was that $P_1$ was the projection that included the product projection $P_C$ or, in the case of the rank-1 method, the projection to rank-1 summands. In the simple setting the structure projections are then assigned to $P_2$, while in the compound setting these also are also assigned to $P_1$ and $P_2$ is tasked with linear compatibility and orthogonality among matrices.

There has been much study of iterative algorithms built from two projections for problems where both of the corresponding constraint sets are convex. Since we will be interested in applications where at least one of the constraint sets is nonconvex, we are limited to schemes that have proven successful even in that setting. One of these is the \textit{alternating direction method of multipliers} or ADMM iteration \cite{B}:
\begin{eqnarray}
x_1&=&P_1(x_2+x)\\
x_2'&=&P_2(x_1-x)\nonumber\\
x'&=&x+\alpha(x_2'-x_1).\nonumber
\end{eqnarray}
Three sets of the original variables are updated in each iteration: $x$, $x_1$ and $x_2$. If in one iteration it happens that $x_1=x_2'$, then $x$ is unchanged and neither are $x_1$ and $x_2$ in the next round. Since $x_1=x_2'=x^*$ is fixed by both projections, we see that ADMM finds a solution whenever it arrives at a fixed point.

By means of the $x$ variables and the positive parameter $\alpha$, the ADMM algorithm is able to escape the traps that plague the more naive algorithm, where the two projections are simply alternated. The traps in the latter algorithm, which is also the $\alpha\to 0$ limit of ADMM (with initialization $x=0$), correspond to pairs of distinct, proximal points $(x_1^*,x_2^*)$ on the two constraint sets. In the presence of such a trap, $x$ is incremented by $\alpha(x_2^*-x_1^*)$ in each iteration and, for $\alpha>0$ and enough iterations, can liberate the algorithm from the trap by re-centering the two projections. The third line of the ADMM update shows that $x$ acts like an accumulator for the discrepancy between constraints.

In this study we will be using a different scheme called \textit{relaxed-reflect-reflect} or RRR \cite{BCL, ABT1, ABT2, E2}. This too is best displayed as an update rule for three sets of variables:
\begin{eqnarray}\label{RRR}
x_1&=&P_1(x)\\
x_2&=&P_2(2 x_1 - x)\nonumber\\
x'&=&x+\beta(x_2-x_1).\nonumber
\end{eqnarray}
RRR derives its name from the fact that it can be compactly written as a relaxed combination of $x$ and constraint-reflections,
\begin{equation}
x'=(1-\beta/2)x+(\beta/2)R_2\circ R_1(x),
\end{equation}
where
\begin{equation}
R_i(x)=2P_i(x)-x,\qquad i\in\{1,2\}.
\end{equation}
With suitable definitions of variables, ADMM with $\alpha=1$ can be shown to be equivalent to RRR with $\beta=1$. The fixed-point/solution relationship for RRR is exactly as it is for ADMM, as are some other features. A relatively minor difference is the fact that for ADMM one must initialize $x$ and $x_2$, compared to just $x$ for RRR. This is truly insignificant for the intended applications, where the variables enter into a rather chaotic steady state dynamics very quickly, thereby losing all memory of the initial conditions. For ADMM it is common practice to initially set the `accumulated discrepancy' variables $x$ to zero.

Once the iteration scheme is selected, there are two ways to optimize the algorithm. While local fixed-point convergence holds for a wide range of the parameters $\alpha$ and $\beta$ ($0<\beta<2$ for RRR), particular settings may prove advantageous for minimizing the much longer times the algorithm spends searching, chaotically, for the fixed-point's basin. A common strategy in global optimization is to combine rounds of different methods, or a schedule of random restarts. Such strategies will have little effect on ADMM/RRR precisely because of the strongly mixing character of the dynamics. Finally, one should consider swapping $1\leftrightarrow 2$ in the ADMM/RRR update rules, as that gives an inequivalent algorithm. 

Our reporting of the RRR algorithm on a sampling of matrix decomposition problems will mostly feature the time series of the root-mean-square constraint discrepancy defined as
\begin{equation}
\Delta=\frac{1}{\sqrt{M}}\,\|x_1-x_2\|_2,
\end{equation}
where normalizing by the number of variables $M$ makes it easier to compare problem instances differing just by size. On hard problems $\Delta$ fluctuates randomly until the variables arrive by chance at the basin of a fixed point, whereupon $\Delta$ decays exponentially to the computer's working precision. The floating point nature of the algorithm usually does not pose a problem, either because the errors in real-world constraint matrices is larger than working precision, or because solutions can be verified with integer arithmetic when there are discrete (\textit{e.g.} integer) structure constraints.

\begin{figure}[t]
\begin{center}
\includegraphics[width=5in]{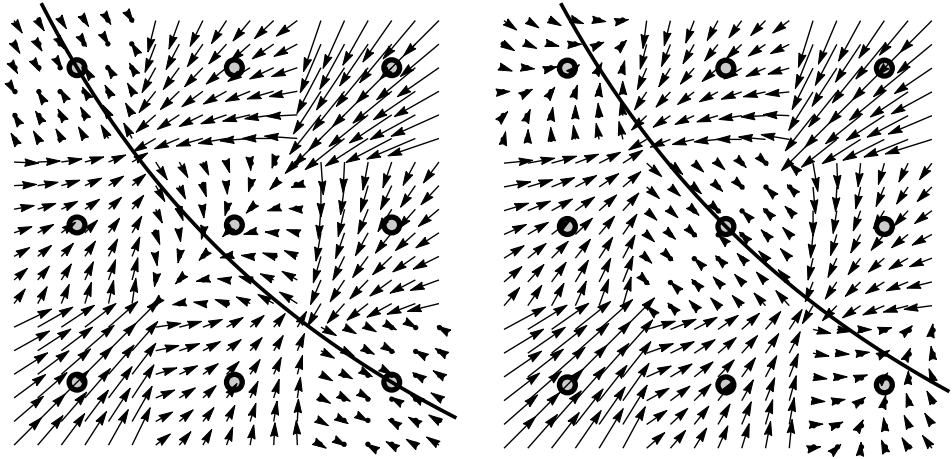}
\end{center}
\caption{Details of the RRR algorithm flow field in a space of two dimensions, where one constraint set is the integer lattice and the other is the curve $x y=c$; \textit{left panel}: $c=15$, \textit{right panel}: $c=16$. In the flow field for $c=15$ there are curves of fixed points passing through the solutions $(3,5)$, $(5,3)$ and limit cycles associated with the near solution $(4,4)$. Fixed points and limit cycles are interchanged in the $c=16$ flow field.}
\end{figure}

The ADMM and RRR algorithm have one potential failure mode when constraint sets are non-convex: rather than converge on fixed points they can get trapped on limit cycles \cite{ABT2}. To better understand the nature of this phenomenon and why it seldom arises in practice, we examine what is probably the first product-constraint problem that comes to mind: integer factorization. The most direct constraint formulation uses the plane $\mathbb{R}^2$ for the factors $(x,y)$, the curve $x y=c$ as one of the constraint sets and $\mathbb{Z}^2$ as the other. We have already seen (Figure 1) how to project to the product constraint, while rounding projects to the integer lattice. To study the dynamics in the plane we examine the flow fields associated with the $\alpha\to 0$, $\beta\to 0$ limits of the update rules. The flow field for the RRR algorithm is the vector field
\begin{equation}
P_1\left(2P_2(x,y)-(x,y)\right)-P_2(x,y),
\end{equation}
and is rendered in Figure 2 for the case that $P_1$ projects to the hyperbola and $P_2$ rounds to the integer lattice. Comparing the flow fields for $c=15$ and $c=16$ (left and right panels of the Figure), we see that the mostly small changes have the effect of transferring curves of fixed points from the $(3,5)$ and $(5,3)$ solutions in one case, to the $(4,4)$ solution in the other. We also see that the fixed point flow near a true solution transforms to a flow field with limit cycles when, by changing $c$, true solutions become near solutions. When we contemplate trying to factor large integers in this constraint formulation, we see that solutions are not very robust to `noise' in the low order bits of the constant $c$, and limit cycles are an unfortunate by-product of this sensitivity.

Arag\'on Artacho and coworkers \cite{ABT2} give other instances of RRR limit cycle pathologies, also in the plane. A reasonable hypothesis that would explain why the phenomenon is not prevalent in applications is the fact that usually many dimensions are required to formulate a problem in terms of constraints, and consequently relatively few bits of information are imposed per dimension. In such formulations the integrity of solutions is robust to noise and there is no need to have many limit cycles that can easily be transformed to fixed points, depending on the vagaries of the noise. Though lacking theoretical support for this hypothesis, we should be wary of applying ADMM or RRR in situations (\textit{e.g.} integer factorization by constraints in the plane) that require high precision in any coordinate of the constraint embedding.

\section{A sampling of applications}

The purpose of this section is to survey the broad range of applications made possible by matrix product constraint projections. By separating the product constraint from structural constraints, projection methods provide a degree of flexibility absent in many other methods. Although it will be clear that projection methods are very efficient for some of the applications, this survey falls short of a comprehensive comparison with alternative methods.

\subsection{Gram matrix decomposition}

In the \textit{maximum determinant problem} one seeks matrices $X\in \{-1,1\}^{m\times m}$ that achieve the highest possible determinant. One strategy for finding such $X$ is to first limit the possible Gram matrices $C=X X^T$ that a maximum determinant $X$ could have. For example, one of the four candidate Gram matrices obtained for $m=15$ \cite{O} had the form $C=12 I_{15}+B$, where $B$ is obtained by removing the last row and column of the matrix
\begin{equation}
\left[
\begin{array}{cccc}
3J&-J&-J&-J\\
-J&3J&-J&-J\\
-J&-J&3J&-J\\
-J&-J&-J&3J
\end{array}
\right],
\end{equation}
where $J$ represents a $4\times 4$ block of 1's. We can try to obtain $X$ from $C$, if it exists, by using the symmetric product constraint projection \eqref{symProj} for one of the two projections in the RRR scheme, and element-wise rounding to $\pm 1$ for the other.

As a warm-up, especially given the uncertainty in the existence of the decomposition, we can construct soluble $m=15$ instances by forming Gram matrices from random $X$ whose elements are uniformly sampled from $\{-1,1\}$. We will use $P_1$ for the discrete structure of the factors and $P_2$ for the smooth space of orthogonal matrices that parameterize the product constraint. This assignment of discrete/smooth constraints in RRR and $\beta=0.2$ worked well on the \textit{bit retrieval problem} \cite{E2}, a special case of Gram matrix decomposition where the matrices have a circulant structure.

\begin{figure}[t]
\begin{center}
\includegraphics[width=5in]{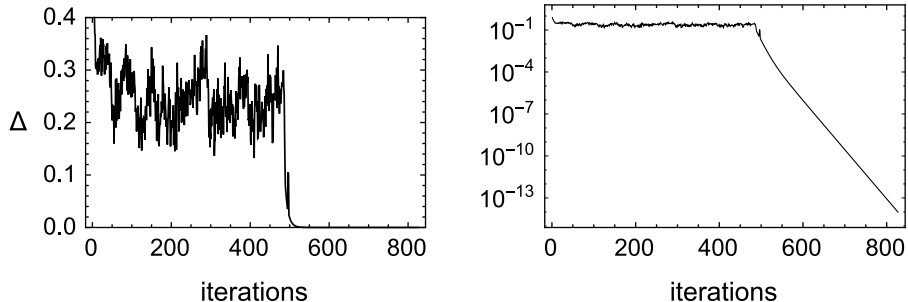}
\end{center}
\caption{Constraint discrepancy time series (log-scale in right panel) in a random $15\times 15$ instance of reconstructing a $\pm 1$ matrix $X$ from its Gram matrix $X X^T$.}
\end{figure}

Not surprisingly, especially given the relationship to bit retrieval, we find there is a strong relationship between the Gram matrix determinant and the number of RRR iterations we should expect before a solution is found. Our random, soluble instances have of course much smaller determinant than the candidate Gram matrices for the maximum determinant problem. The RRR constraint discrepancy time series for a typical one is shown in Figure 3. There is an abrupt change in behavior from `chaotic search', in the first few hundred iterations, to `systematic refinement', in the final iterations. Because the constraint sets in the refinement phase are well approximated by convex sets, the linear convergence we see in the log-discrepancy plot is exactly what we expect of an algorithm designed for convex problems. More remarkable is the fact that the algorithm continues to be reliable, in a statistical sense, even for highly non-convex constraint sets such as we have here. While we do not know when the algorithm will stumble into the attractive basin of a solution and start refining, the statistics of these events have the simplicity of radioactive decay.

Extensive experiments with bit retrieval \cite{E2} show the RRR run times (iteration counts) on fixed instances with random initial $x$ have an exponential distribution. Our (successful) experiments decomposing the proposed maximum determinant Gram matrix, though more limited, are consistent with this property. All 20 attempts produced solutions; the mean iteration count was $5\times 10^5$.

\subsection{Factoring cyclic polynomials}

The problem of factoring polynomials with integer coefficients into polynomials of the same kind, for which there are efficient algorithms \cite{LLL}, is made much harder when posed in the ring of \textit{cyclic polynomials}. The latter is the quotient ring $Z_m=\mathbb{Z}(q)/(q^m-1)$, where exponents are equivalent modulo $m$, for some integer $m$. For example, the polynomial
\begin{equation}
1+2 q^3+3 q^4
\end{equation}
is irreducible in $\mathbb{Z}(q)$ but factors as
\begin{equation}
(1+q+q^4)(1-q+q^3+q^4)
\end{equation}
in the ring $Z_5$. The security of cryptographic keys in protocols such as NTRU \cite{GS} rests in part on the hardness of factoring in $Z_m$.

The problem of factoring a polynomial $c(q)=x(q)y(q)$ in $Z_m$ is equivalent to factoring an $m\times m$ circulant matrix $C$ into circulant matrices $X$ and $Y$. The top rows of the matrices are the polynomial coefficients,
\begin{equation}
c(q)=\sum_{k=0}^{m-1} c_k\, q^{k}\qquad C_{i j}=c_{(j-i\mod{m})},
\end{equation}
and similarly for $x(q)$ and $y(q)$.
By far the most direct way to express the matrix product constraint for circulant matrices is in terms of the Fourier transforms of the polynomial coefficients. Defining these as
\begin{equation}
\hat{c}_l=\frac{1}{\sqrt{m}}\sum_{l=0}^{m-1} e^{i 2\pi k l/m} c_k,\qquad 0\le l\le m-1,
\end{equation}
and similarly for $x$ and $y$, the constraint $X Y=C$ takes the form
\begin{equation}\label{fourierprod}
\hat{x}_l\, \hat{y}_l=\hat{c}_l,\qquad 0\le l\le m-1.
\end{equation}
We recognize this as $m$ independent instances of the complex-scalar product constraint for which the projection was worked out in section \ref{subsubsec:scalar}. That projection, when extended to $m$ independent scalar pairs, minimized the squared distance
\begin{equation}
\sum_{l=0}^{m-1} |\Delta\hat{x}_l|^2+|\Delta\hat{y}_l|^2,
\end{equation}
which equals the squared distance we are using for our circulant matrices,
\begin{equation}
\sum_{k=0}^{m-1} |\Delta x_k|^2+|\Delta y_k|^2,
\end{equation}
by a Fourier transform identity. Note that since the polynomial coefficients are real, the Fourier transforms come in complex-conjugate pairs ($\hat{c}_l$ and $\hat{c}_{-l}$), thereby reducing the number of projections by a factor of two.

To factor polynomials in $Z_m$ by projections, we first embed our polynomials in the ring $R_m=\mathbb{R}(q)/(q^m-1)$. The projection to elements of $Z_m$ is accomplished by rounding all coefficients to the nearest integer. The other projection restores the product constraint by a sequence of three steps: (1) Fourier transforming the coefficients of $x(q)$ and $y(q)$, (2) performing $m$ projections on pairs of Fourier coefficients to the complex-scalar product constraint \eqref{fourierprod}, and (3) inverse Fourier transforming the projected Fourier coefficients to produce a pair of polynomials in $R_m$ that satisfy $x'(q)y'(q)=c(q)$.

As an interesting test of cyclic polynomial factoring by projections, we restrict the coefficients of $x(q)$ and $y(q)$ to be $\pm 1$. For these instances we have a simple upper bound of $2^m$ on the complexity, since by exhausting on the coefficients of $x(q)$, the coefficients of $y(q)$ are found by solving linear equations and then checking for membership in $\{-1,1\}$. Also, we believe the most interesting case is factoring $c(q)$ with small coefficients. The product we will use in our experiments is the $m=23$ polynomial:
\begin{eqnarray}\label{c23}
c(q)&=&1-3q -3q^2 -3q^3 +q^4 +q^5 +q^6 +q^7 +q^8\\
&& -3q^9 -3q^{10} -3q^{11} +q^{12}-3q^{13} -3q^{14}+q^{15}\nonumber\\
&& -3q^{16} -3q^{17} +q^{18} +q^{19} -3q^{20} +q^{21} +q^{22}.\nonumber
\end{eqnarray}
Because the coefficients of the factors are $\pm 1$, all the coefficients of $c(q)$ must be odd and in the same residue class mod 4. The coefficients of non-trivial $c(q)$ that are as small as possible will therefore be two-valued, in this case $-1$ or $3$.

Products $c(q)$ with small coefficients are interesting because they go furthest in probing the non-compact nature of the product constraint. Consider the Fourier-power vectors of the factors: $f_l=|\hat{x}_l|^2$, $g_l=|\hat{y}_l|^2$. Since $\sum_l f_l=\sum_l g_l=m$, these lie in a simplex with vertices on the axes of the positive orthant. When all the coefficients of $c(q)$ are as small as possible, the same holds true of its Fourier coefficients and in particular, the total Fourier power $\sum_l |\hat{c}_l|^2$ is minimized. Since the latter is the inner product $\sum_l f_l g_l$, by minimizing the Fourier power in $c(q)$ we force the power vectors $f_l$ and $g_l$ to have a large separation on the simplex. In terms that matter to the projection algorithm, a large simplex separation translates to many pairs $(\hat{x}_l,\hat{y}_l)$ in the solution with very different magnitudes, \textit{i.e.} points in the `asymptotes' of the constraint `hyperbola'.

To factor \eqref{c23} we used the RRR algorithm with update rule \eqref{RRR}, where $P_1$ is the product constraint projection and $P_2$ projects the polynomial coefficients to $\pm 1$. A factorization was obtained on all attempts with $\beta=0.2$, the same $\beta$ that does well on bit retrieval \cite{E2}. Bit retrieval corresponds to the case of symmetrical factors, $y(q)=x(1/q)$, where projection to the product constraint is the elementary map \eqref{phaseProj} that takes a complex number to the nearest point on a circle. In the non-symmetrical problem the projection is computed by iterating $T$ cycles of quasiprojections and tangent-space projections (section \ref{subsubsec:scalar}). By increasing $T$ we improve the quality of the projection. While increasing $T$ certainly helps fixed-point convergence in the final stage of the solution process, the benefits of a high quality projection in the long, chaotic fixed-point search is less obvious.

With $T=0$, where tangent-space refinement of the projection is turned off, the mean iteration count over 20 trials was 49,000. Adding one cycle of refinement ($T=1$) reduces this to 21,000. Beyond this (20,000 mean iterations for $T=2$) the improvement does not make up for the extra work in computing the projection. We will see that the $T$-dependence of results is more pronounced in other applications. Our results for the mean number of iterations look encouraging relative to the complexity upper bound given by $2^{23}$ linear equation problems.

\subsection{Non-negative matrix factorization}

Applications of non-negative matrix factorization range from small, handcrafted problems in computational geometry and communication complexity, to large-scale industrial problems in data mining and machine learning. In the latter applications an exact factorization usually does not exist, and the task is to find the best approximate factorization. Projection methods, with little modification over how they are used for exact factorization, can also be used in this context. Rather than finding a true fixed point, when there is no exact factorization the ADMM and RRR algorithms are good at finding pairs of proximal points on the two constraint sets \cite{BCL}. One of these points corresponds to matrices with only non-negative entries, and its proximity to the other set implies that the product constraint is nearly satisfied.

In large scale applications the distinction between rank-limited and rank-exces\-sive factors does not come up. In fact, usually the opposite is true: the rank of the approximate factors is required to be smaller, by choice of the middle dimension, than the rank of the (noisy) constraint matrix. Another significant consideration for large scale applications is the fact that the matrices are usually too large to be manipulated as actual matrices. A very different mode of computation, called online learning, is required for these problem.

For the reasons just described, the non-negative matrix factorization problems we consider are of the exact and small variety, as in the recent study by Vandaele \textit{et al.} \cite{VGGT}. The existence of hard problems in this domain became clear when Vavasis \cite{V} showed that determining the non-negative rank of a non-negative matrix is NP-complete. For a non-negative matrix $C\in \mathbb{R}^{m\times n}$ to have non-negative rank $r_+$, it must be possible to express it as the product of a non-negative $X \in \mathbb{R}^{m\times r_+}$ and non-negative $Y\in \mathbb{R}^{r_+\times n}$. We will consider two problems. In the first, $r_+$ is known to equal the standard or real-rank of $C$ and the rank-limited compound projection method of section \ref{sec:ranklim} can be used. The second application features the linear Euclidean distance matrices already introduced in section \ref{sec:exrank}, where the rank-excessive method is required. The latter will be compared with the rank-1 method (section \ref{sec:rank1}) which places no restrictions on the factors.

\subsubsection{Designed instances with zero elements}\label{sec:designed}

For testing algorithms one can generate exact non-negative matrix factorization instances by (1) selecting the matrix shapes $m=n>k$, (2) generating the matrix entries of a solution $(X,Y)$ by uniformly sampling the interval $[0,1]$, and (3) computing the constraint matrix $C=X Y$. However, such instances are easy and do not rigorously test algorithms. We will generate significantly harder instances by forcing a particular fraction of the entries in $X$ and $Y$ to be exactly zero.

To determine the fraction of zeros in $X$ and $Y$ that gives hard instances, we consider the size of the space of solutions. For any instance the space of solutions always contains orbits under the group $G$ of $k\times k$ matrices generated by all permutation matrices as well as arbitrary positive diagonal matrices. This group comprises only non-negative matrices, and for any $g\in G$, the transformed matrices $X g$ and $g^{-1}Y$ give another non-negative factorization. Easy instances are characterized by not just having a single $G$-orbit of solutions, but a continuous space of distinct orbits.

To probe the space of solution orbits we consider the $k\times k$ matrices infinitesimally close to the identity that generate them. Starting with the factorization $C=X Y$, consider the factorization $C=X'Y'+O(\epsilon^2)$ where $X'=X(I_k+\epsilon A)$, $Y'=(I_k-\epsilon A)Y$, and $A$ is an arbitrary $k\times k$ matrix. When $X$ and $Y$ have no zeros, then for small enough $\epsilon$ neither will $X'$ and $Y'$. The space of solutions in that case has $k^2$ generators. Now suppose that a fraction $f$ of the entries in $X$ and $Y$ are zero. The condition that $X'$ remain non-negative translates to a set of linear homogeneous inequalities $(XA)_{i j}\ge 0$, one for each $(i,j)$ where $X_{ij}=0$. Combined with the analogous inequalities that apply to $Y'$, there are in total $M=2 f k m$ inequalities on the $N=k^2$ entries of $A$. In the limit of large matrices, where it is not unreasonable to model the directions that define these inequalities as uniform on the $(N-1)$-sphere, there is a sharp transition\footnote{This is equivalent to the behavior of the probability that $M$ random points on the $(N-1)$-sphere all lie within the same hemisphere, an old problem apparently first analyzed by Ludwig Schl\"afli.} from a cone of feasible $A$, to just $A=0$, when $M/N= 2$. Taking a cue from hardness transitions in other problems \cite{HHW}, we use this onset of uniqueness, where the space of solutions collapses to a single $G$-orbit, as the signal for the hardest kind of instance. This gives $f=k/m$ as the zero fraction for hard problems.

We now present some results for a single random instance of the type described above with $m=n=50$, $k=25$, and $f=1/2$ for the zero fraction. After generating $X$ and $Y$, the product $C=X Y$ was checked to have rank 25. We used the compound projection scheme of section \ref{sec:ranklim}, with $P_1$ in the RRR algorithm combining non-negativity projection on $X$ and $Y$ with the product constraint projection on the $25\times 25$ matrices $W$ and $Z$. The other projection, $P_2$, restores the linear constraints \eqref{lin} that involve the matrices $U$ and $V$ from the singular value decomposition of $C$. As these introduce the metric scale parameters $g$ and $h$, one of our first objectives is to study how the algorithm is affected by them. We keep $g=h$ because our two factors have the same shape.

\begin{figure}[t!]
\begin{center}
\includegraphics[width=3.5in]{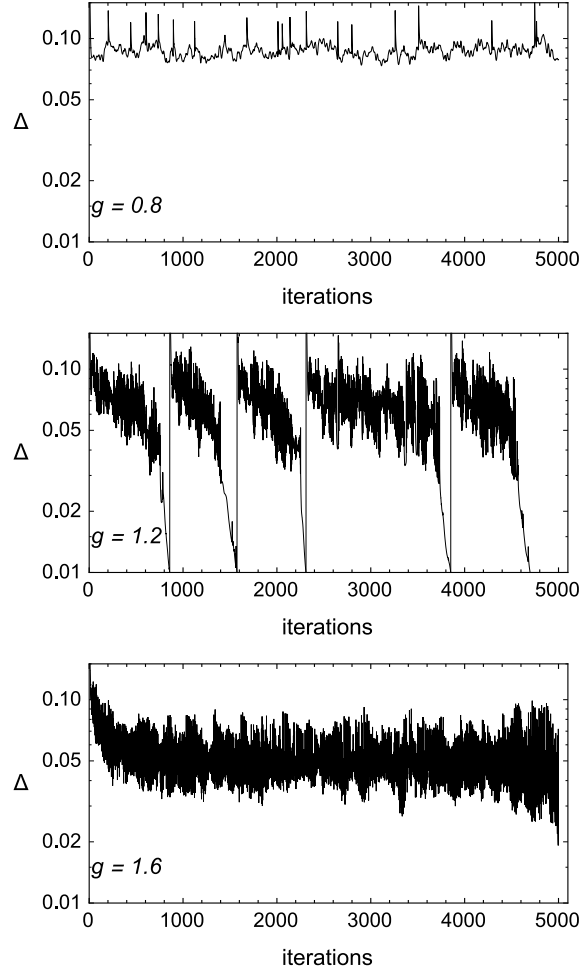}
\end{center}
\caption{Constraint discrepancy time series for a designed instance of non-negative matrix factorization for three values of the metric parameter $g$. Non-negativity is given greater weight than the product constraint when $g$ is small (top panel), and the reverse holds when $g$ is large (bottom panel). The best setting of $g$ is when neither constraint dominates (middle panel, five solutions).}
\end{figure}

With $\beta=0.2$ and the number of tangent-space refinement cycles set conservatively at the high value $T=10$ (see below), the behavior of the RRR constraint discrepancy upon changing the metric parameter $g$ is shown in Figure 4. Not surprisingly, performance degrades both when $g$ is too small and too large. At the optimal value $g\approx 1.2$ the compatibility between $X$ and $W$ (respectively $Y$ and $Z$) is not dominated by one or the other, that is, non-negativity and the product constraint have comparable roles in the search for the solution. All trials with $g=1.2$ produced solutions. A steady, fluctuating behavior of $\Delta$ followed by a sudden drop is characteristic of combinatorial searches when the solution is unique or nearly unique. The factors found by the algorithm (after normalizing columns/rows) proved to be (column/row) permutations of the factors used to create the problem instance.

\begin{figure}[t!]
\begin{center}
\includegraphics[width=3.5in]{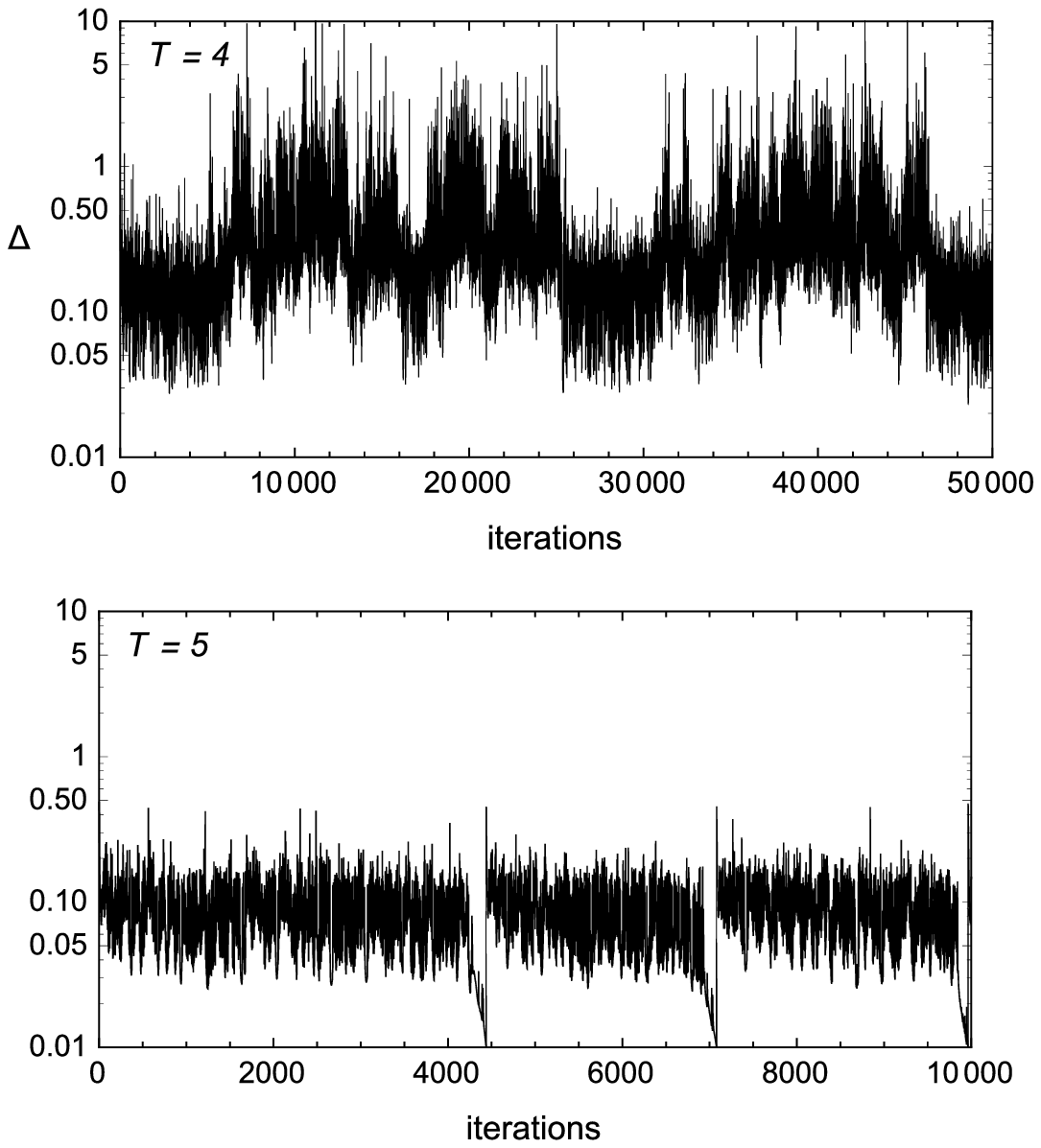}
\end{center}
\caption{Change in the behavior of the constraint discrepancy time series, in a designed instance of non-negative matrix factorization, between $T=4$ cycles of tangent-space refinement and $T=5$. Solutions are found consistently within about 2,000 iterations for $T=5$ (bottom panel) but essentially never when $T=4$.}
\end{figure}

Non-negative matrix factorization makes somewhat higher demands on tangent-space refinement of the constraint projection than what was needed for the scalar products in the cyclic polynomial factorization problem. Fixing $g=1.2$ on the same instance studied above, Figure 5 shows the rather abrupt change in behavior of the discrepancy time series between algorithms with $T=4$ and $T=5$ cycles. With only 4 cycles of refinement the algorithm failed to find a solution in 50,000 iterations, even while showing no sign of getting trapped. We interpret this as a sign that the distance-minimizing quality of the product constraint projection is so poor at $T=4$ that the attractive basins of the RRR fixed points are so small that they have become needles in a haystack. But already with $T=5$ the algorithm consistently finds factorizations, with mean iteration count 2,100. By $T=10$ the mean iteration count is 1,000 and remains at essentially this value for higher $T$. This shows that a critical number of tangent-space refinements of the product constraint projection are essential for the algorithm to work, but that increasing this number beyond that threshold brings diminishing returns.

Vandaele and coworkers \cite{VGGT} proposed a very different family of matrices for testing algorithms, inspired by a problem in communication complexity. These are designed to have the same sparsity pattern as \textit{unique disjointness matrices} and have factors with the following block-substitution rules:
\begin{equation}
X_{d+1}=\left[
\begin{array}{ccc}
X_d&X_d&X_d\\
0&X_d&0\\
X_d&0&0\\
0&0&X_d\\
\end{array}
\right]
\quad
Y_{d+1}=\left[
\begin{array}{cccc}
Y_d&Y_d&0&0\\
Y_d&0&Y_d&0\\
Y_d&0&0&Y_d\\
\end{array}
\right].
\end{equation}
With $X_1=Y_1=I_1$, we see that the instance with constraint $C_d=X_d Y_d$ has factors with shapes $m=n=4^{d-1}$, $k=3^{d-1}$. By inspection we can verify that the factors have equal real and non-negative ranks, and that these match the middle dimension $k$.

\begin{figure}[t]
\begin{center}
\includegraphics[width=4in]{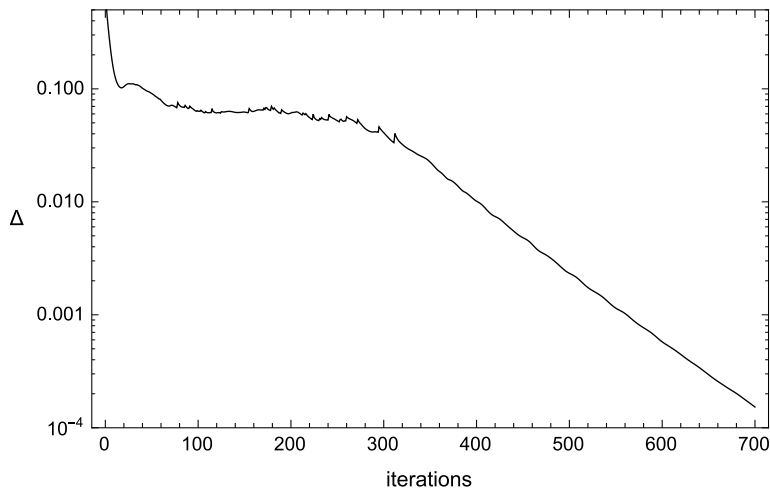}
\end{center}
\caption{Convergent behavior of the RRR discrepancy in an easy case of non-negative factorization based on unique disjointness matrices \cite{VGGT}.}
\end{figure}

The matrices $C_d$ do not pose much of a challenge to the rank-limited compound projection method. With settings $\beta=0.2$, $g=h=0.8$ and $T=10$, the RRR discrepancy for $C_5$, shown in Figure 6, is nearly monotonic-decreasing already in the earliest iterations. The direct passage to convergent behavior is probably a direct consequence of the strong hierarchy of the singular values of these matrices.

\subsubsection{Linear Euclidean distance matrices}

 The linear Euclidean distance matrix $C_m$ of order $m$ has elements
 \begin{equation}
 (C_m)_{i j}=(i-j)^2,\qquad 1\le i\le m,\; 1\le j\le m.
 \end{equation}
These matrices have (for $m\ge 3$) rank $3$ and logarithmically growing non-negative rank $r_+$ \cite{H}. An upper bound $k$ on $r_+$ is given by the middle dimension in a non-negative factorization of $C_m$. As we have no reason to believe the ranks of the factors equal $3$, the rank-excessive construction must be used. In this method one part of each factor, $X_C$ and $Y_C$, has $\rank{(C_m)}=3$. The other part, $X_\perp$ and $Y_\perp$, increases the rank of the factors and is subject to an orthogonality constraint. The two parts of each factor are required to be non-negative when summed and in general are not non-negative individually.

The RRR algorithm can run afoul of limit cycle behavior in this application. With $\beta=1$ (the mid-point of the nominal range) and metric parameters $g=h=0.5$ --- settings that often and quickly lead to solutions --- the algorithm occasionally finds itself in quasi-limit cycles. Although these are unstable and do not represent permanent traps, the search performed by the algorithm during these epochs is not very productive. An example from an attempted $k=5$ factorization of $C_6$ is shown in Figure 7.

A tendency for limit cycles is consistent with the general caution of section \ref{sec:iterproj}, that the constraints to combinatorially hard problems should not require a large number of bits of information per Cartesian dimension of the constraint-space. Here the principle would apply to the $k$ rank-1 summands $Z^l_{i j}=X_{i l}Y_{l j}$ (fixed $l$) whose sum must give a partition of the integers in $C_m$, ranging from $0$ to $(m-1)^2$, into integers.

\begin{figure}[t]
\begin{center}
\includegraphics[width=4in]{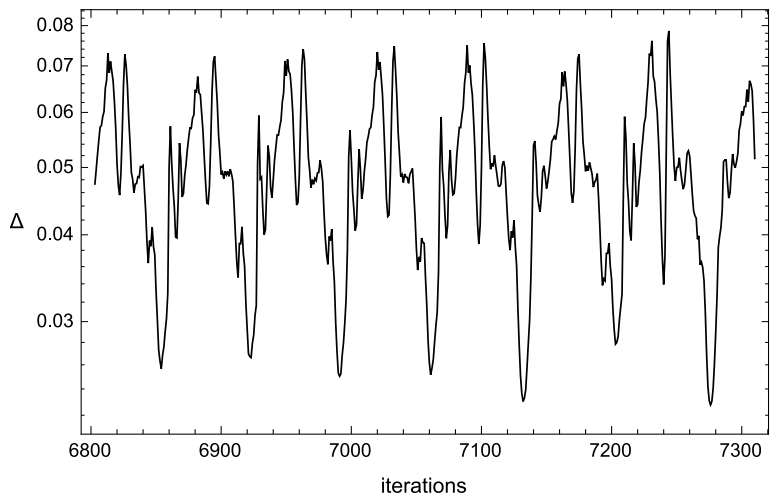}
\end{center}
\caption{Quasi-limit cycle behavior in a non-negative factorization of the order 6 linear Euclidean distance matrix.}
\end{figure}

Through experimentation we found that limit cycle behavior can be avoided by using a reasonable initial point for the RRR algorithm. Define the SVD-based factorization as $X_C=U\sqrt{D(r,k)}$, $Y_C=\sqrt{D(k,r)}\,V$, where the diagonal matrix of singular values $D$ has been extended with zeroes to have the correct shape. To produce non-negative factors, define $X_\perp=\max{(0,-X_C)}$, $Y_\perp=\max{(0,-Y_C)}$. The  point (in the rank-excessive construction)
\begin{equation}
(\sqrt{D(r,k)},X_C,X_C,X_\perp,X_\perp\;;\;\sqrt{D(k,r)},Y_C,Y_C,Y_\perp,Y_\perp)
\end{equation}
satisfies all constraints except the orthogonality property \eqref{ortho3}. Running RRR with this as initial point and the parameters above, a non-negative $k=5$ factorization of $C_6$ is found in 786 iterations. For these factorizations we terminate the algorithm when the summand matrices $Z^l$, after rounding to integer matrices, are rank-1 and sum to $C$. The $k=6$ factorization of $C_8$ required 1,508 iterations and $k=7$ for $C_{12}$ required 88,467. For $C_{16}$ the search was found to be more productive with metric parameters $g=h=0.3$. The discrepancy time series of a successful $k=8$ factorization in 53,007 iterations is shown in Figure 8. In all of these experiments the value of $k=r_+$ is the smallest possible. The ranks of the factors in this sequence of instances grows as $(4,4)$, $(4,5)$, $(5,5)$, $(5,6)$.

\begin{figure}[t]
\begin{center}
\includegraphics[width=4in]{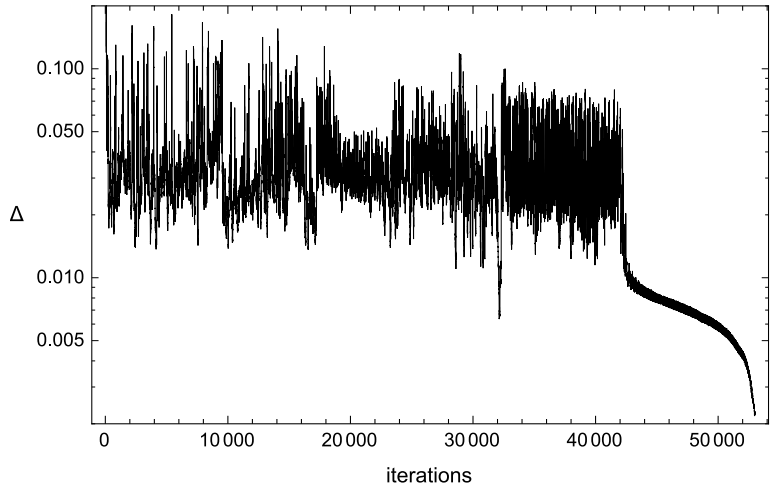}
\end{center}
\caption{RRR discrepancy in a successful non-negative factorization of the order 16 linear Euclidean distance matrix.}
\end{figure}

The rank-1 method may also be used for these instances of non-negative factorization. As this method works in a space with more dimensions than the product-constraint method (for large $m$), there is reason to hope the limit cycle problem will be mitigated. This turns out not to be the case. Using the projections \eqref{rank1P1} and \eqref{rank1P2} in the RRR algorithm with $\beta=1$ on the $k=5$ factorization of $C_6$, we observe trapping on limit cycles that appears to be permanent in about 20\% of trials. In these trials the initial random $Z^l_{i j}$ elements are uniform samples between $0$ and $(m-1)^2$. The mean iteration count in the untrapped trials is 1,600, about twice the number needed by the product-constraint method. It appears the limit cycle problem is mitigated by replacing the simplex projection $P_*$ for the structure by the stronger projection $P_A\circ P_*$ that makes use of the fact that in these instances the $Z^l_{i j}$ are integers. For the $k=r_+$ factorizations of $C_6$ and $C_8$ the algorithm now averages 560 and 5,000 iterations; $C_{12}$ and $C_{16}$ are still out of reach.

\subsubsection{Comparison with norm minimization methods}

State-of-the-art non-negative factorization methods \cite{VGGT} are all based on the minimization of
\begin{equation}\label{objective}
\|X Y-C\|_2,
\end{equation}
differing only on strategies for solving this non-convex optimization problem. The latter  include alternating a sequence of non-negative minimizations with respect to one factor while the other is held fixed, or a similar strategy applied to individual rows/columns of the factors. As these restricted convex minimizations invariably arrive at non-zero local minima of the objective \eqref{objective}, a significant degree of randomization is required for these methods to succeed. The best strategies \cite{VGGT} in that regard involve local randomization, similar in spirit to what is done in simulated annealing. By contrast, the only explicit randomness invoked by projection methods is in the selection of the initial point. But as we have seen, in the case of the linear Euclidean distance matrices even this degree of randomness is unnecessary as a well motivated special initial point achieves good results.

The assertion that projection methods are just another technique for global optimization neglects a number of possibly relevant points. First, we note that non-negative factorization by minimization of \eqref{objective} never makes use of the fact that, in exact problems, the minimum of the objective is zero. This fact plays a central role in developing projection methods for this problem. A second point is that non-negative factorization problems may have interesting structure that minimization methods do not exploit. For example, we are not aware of minimization methods that address the two cases of the (real) ranks of the factors (rank-limited vs. rank-excessive), as we were forced to consider in the construction of compound projections. Lastly, minimization methods normally are unable to take advantage of discrete constraints (integer, 0-1) on the factors (or rank-1 summands).

\section{Summary}

Fast projections to the matrix product constraint enables new methods for finding matrices that not only have a given product but also have a particular structure (\textit{e.g.} non-negativity). The first step in implementing these methods is to determine if the shapes and ranks of the factors are amenable to one of the simple projections (section \ref{sec:simple}) or whether one of the compound constructions involving additional matrices is required (section \ref{sec:compound}). All of these projections are built from standard matrix decomposition algorithms (Cholesky, singular value, eigenvalue). The core algorithm for most of these projections (section \ref{sec:outerfullrank}) iterates a quasiprojection to the true constraint with a true projection to the tangent-space approximation of the constraint to get a high quality projection.

Once the product and structure constraints are implemented as projections, always as a pair comprising simple or compound projections, an iterative projection method such as ADMM or RRR is used to find matrices that are fixed by both projections and therefore solve the problem. Whereas convergence results for these iterative methods is limited to problems with convex constraint sets, their success with non-convex, combinatorially hard problems makes them an attractive heuristic in that domain. This work examined the strengths and weaknesses of these methods in a variety of problems, including Gram matrix decomposition, factoring cyclic polynomials, and non-negative matrix factorization.

We have not carried out systematic benchmarks for comparison with other global optimization methods, but instead have used our selection of applications to highlight features that for the most part are unique to projection methods. Not least of the questions confronting first-time users is the selection of parameters. Probably the most important are the metric parameters. These appear only in the compound setting (section  \ref{sec:compound}) and determine the distance scales that are applied to all the matrices in the construction. We showed in section \ref{sec:designed} that the optimal setting of the $g$ parameter  is such that neither non-negativity nor the product constraint dominates the other.

The refinement cycle number $T$ and RRR parameter $\beta$ are less critical. Our product constraint projections always produce pairs of matrices that have the correct product and fall short of true projections by failing to be distance minimizing. By increasing the number of refinement cycles $T$, the quality of the projections is improved. Our polynomial and non-negative factoring experiments showed that to achieve good results in these combinatorially hard problems it is only necessary for $T$ to exceed a relatively small number. Finally, a recent study of the RRR algorithm in bit retrieval \cite{E2} suggests a similar threshold effect applies to the $\beta$ parameter. The most efficient search performed by RRR appears to be in the regime where the discrete dynamics is approximating the continuous flow of the $\beta\to 0$ limit.

\section{Acknowledgements}\label{sec:ack}

Most of this work was competed on sabbatical at Disney Research, Boston. I thank Disney and the Simons Foundation for financial support during that period.

\end{document}